\documentclass[a4paper,12pt]{amsart}
\usepackage{amsmath,amsthm,amssymb,a4wide}
\usepackage[arrow,matrix]{xy}
\usepackage{hyperref}
\usepackage{ulem}
\usepackage{graphicx}
\usepackage{pstricks,pst-node,pst-tree,pst-plot}

\theoremstyle{ams}
\newtheorem{theorem}{Theorem}[section]
\newtheorem{proposition}[theorem]{Proposition}
\newtheorem{lemma}[theorem]{Lemma}
\newtheorem{corollary}[theorem]{Corollary}

\theoremstyle{definition}

\newtheorem{remark}[theorem]{Remark}

\newtheorem{example}[theorem]{Example}

\newcommand{\R}{\mathbb{R}}
\newcommand{\Z}{\mathbb{Z}}

\begin{document}
\title[The log-concavity conjecture]{The log-concavity conjecture on semifree symplectic $S^1$-manifolds with isolated fixed points}

\author[Y. Cho]{Yunhyung Cho}
\address{School of Mathematics, Korea Institute for Advanced
Study, 87 Hoegiro, Dongdaemun-gu, Seoul, 130-722, Republic of Korea}
\email{yhcho@kias.re.kr}

\date{\today}
\maketitle

\begin{abstract}

      Let $(M,\omega)$ be a closed $2n$-dimensional semifree Hamiltonian $S^1$-manifold with only isolated fixed points. We prove that a density function of the Duistermaat-Heckman measure is log-concave. Moreover, we prove that $(M,\omega)$ and any reduced symplectic form satisfy the Hard Lefschetz property.
\end{abstract}

\section{Introduction}

    Let $T$ be a connected abelian Lie group, i.e. a torus.
    Let $(M,\omega)$ be a closed $2n$-dimensional Hamiltonian $T$-manifold with the moment map $\mu \colon M \rightarrow \textbf{t}^*$, where $\textbf{t}$ is the Lie algebra of $T$. A Liouville measure $m_L$ on $M$ is defined by
    $$m_L(U) := \int_U \frac{\omega^n}{n!} $$ for any open set $U \subset M$.  Then the push-forward measure $m_{\mathrm{DH}} := \mu_* m_L$, called
    the \textit{Duistermaat-Heckman measure}, can be regarded as a measure on $\mathfrak{t^*}$ such that for any Borel subset $B
    \subset \mathfrak{t^*}$, $m_{\mathrm{DH}}(B) = \int_{\mu^{-1}(B)}
    \frac{\omega^n}{n!}$ tells us that how many states of our system have momenta in $B.$
    Due to \cite{DH}, the Duistermaat-Heckman measure is absolutely continuous with respect to the Leabesque measure on $\mathbf{t}^*$ and the corresponding density function, denoted by $\mathrm{DH}$ and called the $\textit{Duistermaat-Heckman function}$, is a polynomial on any regular open set.
    More precisely, for any regular value $\xi \in \textbf{t}^*$, $$\mathrm{DH}(\xi) = \int_{M_{\xi}} \frac{1}{(n-l)!}\omega_{\xi}^{n-l}$$ where $l$ is a dimension of $T$, $M_{\xi}$ is the symplectic reduction at $\xi$, and $\omega_{\xi}$ is the corresponding reduced symplectic form on  $M_{\xi}$.

    In 1996, W.Graham \cite{Gr} proved that if the Hamiltonian $T$-action on a K\"{a}hler manifold $(M,\omega,J)$ is holomorphic, then the Duistermaat-Heckman measure is log-concave, i.e. $\log{\mathrm{DH}}$ is a concave function on the image of the moment map.

    \begin{example}\label{ex1}
    Consider a compact symplectic toric manifold $(M,\omega,T)$ with the moment map $\mu \colon M \rightarrow \textbf{t}^*$. Since any symplectic
    toric manifold can be obtained by the K\"{a}hler reduction, $M$ is a smooth toric variety and the given $T$-action is holomorphic with
    some $\omega$-compatible complex structure on $M$.
    By Atiyah-Guillemin-Sternberg convexity theorem, the image of $\mu$ is a rational convex polytope satisfying the non-singularity condition.
    Choose any circle subgroup $S^1 \subset T$. Then the moment map of the circle action on $M$ is a composition map $\mu_{S^1} \colon M \stackrel{\mu}\longrightarrow \textbf{t}^* \stackrel{\pi}\longrightarrow \textbf{s}^*$ where $\pi$ is a dual map of $i \colon \textbf{s} \rightarrow \textbf{t}$, the inclusion of Lie algebra induced by $S^1 \hookrightarrow T$. The corresponding Duistermaat-Heckman function $\mathrm{DH}(t)$ is just $C \cdot \mathrm{vol}(\mu(M) \cap \pi^{-1}(t))$ with respect to the Lebesque measure on $\textbf{t}^*$, which is log-concave by Graham's theorem \cite{Gr}. Here, $C$ is a global constant depending on the embedding $S^1 \hookrightarrow T$. Therefore, for any given Delzant polytope $\triangle$ and any height function $h(t) \colon \triangle \rightarrow \R$, the corresponding slice volume function $g : h(\triangle) \rightarrow \R$ defined by $g(t) := \mathrm{vol}(h^{-1}(t))$ is log-concave.

     The following figure is the case when $M = S^2 \times S^2 \times S^2$, $\omega = 2\omega_{FS} \oplus 3\omega_{FS} \oplus 6\omega_{FS}$, and a Hamiltonian circle action on $M$ is given by
    $t \cdot (z_1,z_2,z_3) = (tz_1, tz_2, tz_3)$ for any $t \in S^1$.

\begin{figure}[ht]
\centering

        \psset{unit=20pt}

    \begin{pspicture}(0,-9)(27,8)\footnotesize

            \pspolygon[fillstyle=solid, linewidth=1.5pt](2,0)(5,3)(5,4.5)(4,5.5)(1,2.5)(1,1) \psline(9,0)(9,5.5)
            \psline(2,0)(2,1.5) \psline(2,1.5)(5,4.5) \psline(2,1.5)(1,2.5)
            \psdot(4,4)\psdot(2,0)\psdot(1,1)\psdot(5,3)\psdot(5,4.5)\psdot(2,1.5)\psdot(1,2.5)\psdot(4,5.5)
            \psdot(9,0)\psdot(9,1)\psdot(9,1.5)\psdot(9,2.5)\psdot(9,3)\psdot(9,4)\psdot(9,4.5)\psdot(9,5.5)
            \psline[linestyle=dashed,linecolor=black](1,1)(4,4)
            \psline[linestyle=dashed,linecolor=black](4,4)(4,5.5)
            \psline[linestyle=dashed,linecolor=black](4,4)(5,3)
            \rput(2,-0.5){\tiny{$(S,S,S)$}}
            \rput(0.2,1){\tiny{$(N,S,S)$}}
            \rput(3,1.7){\tiny{$(S,N,S)$}}
            \rput(0.2,2.5){\tiny{$(N,N,S)$}}
            \rput(6,3.2){\tiny{$(S,S,N)$}}
            \rput(6,4.7){\tiny{$(S,N,N)$}}
            \rput(3.4,4.2){\tiny{$(N,S,N)$}}
            \rput(4,6){\tiny{$(N,N,N)$}}

            \rput(9.4,0){\tiny{$0$}}
            \rput(9.4,1){\tiny{$2$}}
            \rput(9.4,1.5){\tiny{$3$}}
            \rput(9.4,2.5){\tiny{$5$}}
            \rput(9.4,3){\tiny{$6$}}
            \rput(9.4,4){\tiny{$8$}}
            \rput(9.4,4.5){\tiny{$9$}}
            \rput(9.4,5.5){\tiny{$11$}}
            \rput(12,3.2){\tiny{DH}}
            \rput(7.5,2.4){$\mu$}

            \rput(16.2,-0.2){\tiny{$0$}}
            \rput(17.2,-0.2){\tiny{$2$}}
            \rput(17.8,-0.2){\tiny{$3$}}
            \rput(19,-0.2){\tiny{$5$}}
            \rput(19.6,-0.2){\tiny{$6$}}
            \rput(20.8,-0.2){\tiny{$8$}}
            \rput(21.4,-0.2){\tiny{$9$}}
            \rput(22.6,-0.2){\tiny{$11$}}
            \rput(23.9,-0.2){$\mu(M)$}
            \rput(15,6){$\mathrm{DH}(t)$}

            \psdot(16,0)\psdot(17.2,0)\psdot(17.8,0)\psdot(19,0)\psdot(19.6,0)\psdot(20.8,0)\psdot(21.4,0)\psdot(22.6,0)
            \pscurve(16,0)(16.15,0.03125)(16.3,0.125)(16.45,0.28125)(16.6,0.5)(16.75,0.78125)(16.9,1.125)(17.05,1.53125)(17.2,2)
            \psdot(17.2,2)
            \psline(17.2,2)(17.8,4)
            \psdot(17.8,4)
            \pscurve(17.8,4)(18.1,4.875)(18.4,5.5)(18.7,5.875)(19,6)
            \psdot(19,6)
            \psline(19,6)(19.6,6)
            \psdot(19.6,6)
            \pscurve(19.6,6)(19.9,5.875)(20.2,5.5)(20.5,4.875)(20.8,4)
            \psdot(20.8,4)
            \psline(20.8,4)(21.4,2)
            \psdot(21.4,2)
            \pscurve(21.4,2)(21.7,1.125)(22,0.5)(22.3,0.125)(22.6,0)

            \psline(8,-5)(16,-5)
            \psline(9,-9)(9,-2)
            \rput(7.8,-2){$\log\mathrm{DH}(t)$}
            \rput(17,-5.2){$\mu(M)$}
            \psdot(9,-5)\psdot(10.2,-5)\psdot(10.8,-5)\psdot(12,-5)\psdot(12.6,-5)\psdot(13.8,-5)\psdot(14.4,-5)\psdot(15.6,-5)
            \rput(9.2,-5.2){\tiny{$0$}}
            \rput(10.2,-5.2){\tiny{$2$}}
            \rput(10.8,-5.2){\tiny{$3$}}
            \rput(12,-5.2){\tiny{$5$}}
            \rput(12.6,-5.2){\tiny{$6$}}
            \rput(13.8,-5.2){\tiny{$8$}}
            \rput(14.4,-5.2){\tiny{$9$}}
            \rput(15.6,-5.2){\tiny{$11$}}
            \pscurve(9.15,-8.46573)(9.3,-7.07944)(9.45,-6.26851)(9.6,-5.69314)(9.75,-5.24686)(9.9,-4.88222)(10.05,-4.57392)(10.2,-4.30686)
            \psdot(10.2,-4.30686)
            \pscurve(10.2,-4.30686)(10.5,-3.88)(10.8,-3.6137)
            \psdot(10.8,-3.6137)
            \pscurve(10.8,-3.6137)(11.1,-3.41588)(11.4,-3.29525)(11.7,-3.22929)(12,-3.20824)
            \psdot(12,-3.20824)
            \psline(12,-3.20824)(12.6,-3.20824)
            \psdot(12.6,-3.20824)
            \pscurve(12.6,-3.20824)(12.9,-3.22929)(13.2,-3.29525)(13.5,-3.41588)(13.8,-3.6137)
            \psdot(13.8,-3.6137)
            \pscurve(13.8,-3.6137)(14.1,-3.88)(14.4,-4.30686)
            \psdot(14.4,-4.30686)
            \pscurve(14.4,-4.30686)(14.7,-4.88222)(15,-5.69314)(15.3,-7.07944)(15.45,-8.46573)

            \psline(15,0)(23,0)
            \psline(16,0)(16,6)
            \psline{->}(7,2.8)(8,2.8)
            \psline{->}(11,2.8)(13,2.8)
            \psline[linestyle=dashed,linecolor=red,linewidth=0.2pt](2,0)(9,0)
            \psline[linestyle=dashed,linecolor=red,linewidth=0.2pt](1,1)(9,1)
            \psline[linestyle=dashed,linecolor=red,linewidth=0.2pt](2,1.5)(9,1.5)
            \psline[linestyle=dashed,linecolor=red,linewidth=0.2pt](5,3)(9,3)
            \psline[linestyle=dashed,linecolor=red,linewidth=0.2pt](5,4.5)(9,4.5)
            \psline[linestyle=dashed,linecolor=red,linewidth=0.2pt](1,2.5)(9,2.5)
            \psline[linestyle=dashed,linecolor=red,linewidth=0.2pt](4,5.5)(9,5.5)
            \psline[linestyle=dashed,linecolor=red,linewidth=0.2pt](4,4)(9,4)
    \end{pspicture}
    \caption{Example \ref{ex1}} \label{figure: illustrate}
    \end{figure}
\end{example}

    Note that any holomorphic Hamiltonian $S^1$-manifold $(M,\omega,J)$ with an $S^1$-invariant K\"{a}hler structure $\omega$ satisfies the followings.
    \begin{enumerate}
        \item $\omega$ satisfies the Hard Lefschetz property,
        \item any reduced symplectic form $\omega_t$ is K\"{a}hler, and hence it satisfies the Hard Lefschetz property, and
        \item the Duistermaat-Heckman measure is log-concave.
    \end{enumerate}
    Hence it is natural to ask that when a Hamiltonian $S^1$-manifold $(M,\omega)$ and the symplectic reduction $(M_t, \omega_t)$ satisfy the Hard Lefschetz property, and when the corresponding Duistermaat-Heckman measure is log-concave.

    In this paper, we focus on the case when $(M,\omega)$ is a closed $2n$-dimensional symplectic manifold with a semifree Hamiltonian $S^1$-action whose fixed point set $M^{S^1}$ consists of isolated points. Let $\mu \colon M \rightarrow \R$ be a corresponding moment map which satisfies $i_X\omega = -d\mu$.
    As noted above, the Duistermaat-Heckman function $\mathrm{DH} \colon \mu(M) \rightarrow \R$ is defined by $$\mathrm{DH}(t) = \int_{M_t} \frac{1}{(n-1)!} \omega_t^{n-1}$$ where $M_t$ is the reduced space $\mu^{-1}(t)/S^1$ and $\omega_t$ is the reduced symplectic form on $M_t$.
    Now we state our main results.

    \begin{theorem}\label{first}
         Let $(M,\omega)$ be a closed symplectic manifold with a semifree Hamiltonian $S^1$-action whose fixed point set $M^{S^1}$ consists of isolated points. Then the Duistermaat-Heckman measure is log-concave.
    \end{theorem}

    \begin{theorem}\label{second}
        Let $(M,\omega)$ be a closed semifree Hamiltonian $S^1$-manifold whose fixed points are all isolated, and let $\mu$ be the moment map. Then $\omega$ satisfies the Hard Lefschetz property. Moreover, the reduced symplectic form $\omega_t$ satisfies the Hard Lefschetz property for every regular value $t$.
    \end{theorem}

    In Section 2, we briefly review the Tolman and Weitsman's work \cite{TW1} which is very powerful to analyze the equivariant cohomology of the Hamiltonian $S^1$-manifold with isolated fixed points. Especially we use the Tolman-Weitsman's basis of the equivariant cohomology $H^*_{S^1}(M)$ which is constructed by using the equivariant version of the Morse theory \cite{TW2}.
    In Section 3, we express the Duistermaat-Heckman function explicitly in terms of the integration of some cohomology class on the reduced space. And then we compute the integration by using the Jeffrey-Kirwan residue formula \cite{JK}. Consequently, we will show that the log-concavity of the Duistermaat-Heckman measure is completely determined by
    the set of pairs $\{(\mu(F),m_F)_F | F \in M^{S^1}\}$, where $\mu(F)$ is the image of the moment map of $F$ and $m_F$ is the product of all weights of the $S^1$-representation on $T_F M$. In Sections 4, we will prove Theorem \ref{first}. And we will prove Theorem \ref{second} in Section 5.

    \begin{remark}\label{ree}

        As noted in \cite{O2}, \cite{K} and \cite{Lin}, Ginzburg and Knudsen conjectured that for any closed Hamiltonian $T$-manifold, the corresponding Duistermaat-Heckman measure is log-concave. Note that the log-concavity conjecture is  motivated by a concavity property of the Boltzmann' entropy in statistical mechanics. (See \cite{O3} for the detail.)
        The log-concavity problem of the Duistermaat-Heckman measure is proved by A.Okounkov \cite{O1} when $M$ is a co-adjoint orbit of the classical Lie groups of type $A_n$, $B_n$, or $C_n$ with the maximal torus $T$-action, by using the log-concavity of the multiplicities of the irreducible representation of each Lie groups.
         But the counter-example was found by Y.Karshon \cite{K}. By using the Lerman's symplectic cutting method, she constructed a closed 6-dimensional semifree Hamiltonian $S^1$-manifold with two fixed components such that the Duistermaat-Heckman measure is not log-concave. Later, Y.Lin \cite{Lin} generalized the construction of 6-dimensional Hamiltonian $S^1$-manifolds not satisfying the log-concavity of the Duistermaat-Heckman measure, and also proved the log-concavity conjecture for any Hamiltonian $T^{n-2}$-manifold $(M^{2n},\omega)$ whose reduced space has $b_2^+ = 1$.
    \end{remark}


\section{Tolman-Weitsman basis of the equivariant cohomology $H^*_{S^1}(M;\Z)$}

    In this section, we breifly review Tolman and Weitsman's results in \cite{TW1}.
    Throughout this section, we assume that $(M^{2n},\omega)$ is a closed semifree Hamiltonian $S^1$-manifold whose fixed points are isolated.
    Note that for each fixed point $p \in M^{S^1}$, \textit{the index of $p$} is the Morse index of the moment map at $p$ which is the same as the twice of the number of negative weights of the tangential $S^1$-representation at $p$.

    \begin{proposition}\cite{TW1}
        Let $N_k$ be the number of fixed points of index $2k$. Then $N_k = {n \choose k}$. Hence $N_k$ is the same as the one of the standard diagonal circle action on $(S^2 \times \cdots \times S^2, \omega_1 \oplus \cdots \oplus \omega_n)$, where $\omega_i$ is the Fubini-Study form on $S^2$ of $i$-th factor.
    \end{proposition}

    \begin{theorem}\cite{TW1}\label{basis}
        Let $2^{[n]}$ be the power set of $\{1,\cdots,n\}$. Then there exist a bijection $\phi : M^{S^1} \rightarrow 2^{[n]}$ such that
        \begin{enumerate}
            \item For each index-$2k$ fixed point $x \in M^{S^1}$, $|\phi(x)| = k$.
            \item Let $u$ be the generator of $H^*(BS^1,\Z)$. For each index-$2k$ fixed point $x \in M^{S^1}$, there exists a unique cohomology class $\alpha_x \in H^{2k}_{S^1}(M;\Z)$ such that for any $x' \in M^{S^1}$, \\
                   $\alpha_x|_{x'} = u^k$ if $\phi(x) \subset \phi(x')$. \\
                   $\alpha_x|_{x'} = 0$ otherwise.
        \end{enumerate}
        Moreover $\{\alpha_x | x \in M^{S^1} \}$ forms a basis of $H^*_{S^1}(M,\Z)$.
    \end{theorem}

    Applying Theorem \ref{basis} to $(S^2 \times \cdots \times S^2, \omega_1 \oplus \cdots \oplus \omega_n)$ with the diagonal semifree Hamiltonian circle action, we have a bijection
    $\psi : (S^2 \times \cdots \times S^2)^{S^1} \rightarrow 2^{[n]}$ and there is a basis $\{\beta_y| y \in  (S^2 \times \cdots \times S^2)^{S^1} \}$
    of $H^*_{S^1}(S^2 \times \cdots \times S^2;\Z)$ satisfies the conditions in Theorem \ref{basis}. Hence we have an identification map
    $$\psi^{-1} \circ \phi : M^{S^1} \rightarrow (S^2 \times \cdots \times S^2)^{S^1}$$ and $\psi^{-1} \circ \phi$ preserves the indices of the fixed points.

    Note that $\psi^{-1} \circ \phi$ gives an identification between $H^*_{S^1}(M;\Z)$ and $H^*_{S^1}(S^2 \times \cdots \times S^2;\Z)$ as follow.
    Let $a_i = \alpha_{\phi^{-1}\{i\}} \in H^2_{S^1}(M;\Z)$ and $b_i = \beta_{\psi^{-1}\{i\}} \in H^2_{S^1}(S^2 \times \cdots \times S^2;\Z)$.
    \begin{lemma}\cite{TW1}\label{tol}
        For each $x \in M^{S^1}$, $\alpha_x = \prod_{j \in \phi(x)} a_j$. Similarly, we have $\beta_y = \prod_{j \in \psi(y)} b_j$ for each $y \in
        (S^2 \times \cdots \times S^2)^{S^1}$.
    \end{lemma}
    \begin{proof}
        For an inclusion $i \colon M^{S^1} \hookrightarrow M$, we have a natural ring homomorphism
        $i^* \colon H^*_{S^1}(M) \rightarrow H^*_{S^1}(M^{S^1}) \cong H^*(M^{S^1}) \otimes H^*(BS^1)$. Kirwan injectivity implies that
        $i^*$ is injective. Hence it is enough to show that $\alpha_x|_z = (\prod_{j \in \phi(x)} a_j)|_z$ for all $x,z \in M^{S^1}$.
        Here, the operation $``|_z"$ is the restriction $H^*_{S^1}(M^{S^1}) \rightarrow H^*_{S^1}(z) \cong H^*(BS^1)=\R[u]$ induced by the inclusion
        $z \hookrightarrow M^{S^1}$. For any $x,z \in M^{S^1}$ with $\mathrm{Ind}(x) = 2k$,
        \begin{itemize}
              \item     $\alpha_x|_{z} = u^k$ if $\phi(x) \subset \phi(z)$.
              \item     $\alpha_x|_{z} = 0$ otherwise.
        \end{itemize}
        On the other hand, $(\prod_{j \in \phi(x)} a_j)|_z = \prod_{j \in \phi(x)} a_j|_z$. Since $a_j|_z = u$ if and only if $j \in \phi(z)$, we have
         \begin{itemize}
              \item     $(\prod_{j \in \phi(x)} a_j)|_z = u^k$ if $\phi(x) \subset \phi(z)$.
              \item     $(\prod_{j \in \phi(x)} a_j)|_z = 0$ otherwise.
        \end{itemize}
        Therefore, $\alpha_x = \prod_{j \in \phi(x)} a_j$ by the Kirwan injectivity theorem. The proof of the second statement is similar.
    \end{proof}

     Hence the $H^*(BS^1)$-module isomorphism $f : H^*_{S^1}(M;\Z) \rightarrow H^*_{S^1}(S^2 \times \cdots \times S^2;\Z)$ which sends $\alpha_x$ to $\beta_{\psi^{-1} \circ \phi(x)}$ for each $x \in M^{S^1}$ is in fact a ring isomorphism by the lemma \ref{tol}.
     To sum up, we have the following corollary.


    \begin{corollary}\label{iso}\cite{TW1}
        There is a ring isomorphism $f : H^*_{S^1}(M;\Z) \rightarrow H^*_{S^1}(S^2 \times \cdots \times S^2;\Z)$ which sends $\alpha_x$ to $\beta_{\psi^{-1} \circ \phi(x)}$. Moreover, for any $\alpha \in H^*_{S^1}(M;\Z)$ and any fixed point $x \in M^{S^1}$,
        we have $\alpha|_x = f(\alpha)|_{\psi^{-1} \circ \phi(x)}$.
    \end{corollary}

\section{The Duistermaat-Heckman function and the residue formula}

   Let $(M,\omega)$ be a $2n$-dimensional closed Hamiltonian $S^1$-manifold with the moment map $\mu \colon M \rightarrow \R$. We may assume that $0$ is a regular value of $\mu$ such that $\mu^{-1}(0)$ is non-empty. Choose two consecutive critical values $c_1$ and $c_2$ of $\mu$ so that the open interval $(c_1,c_2)$ consists of regular values of $\mu$ and contains $0$. By the Duistermaat-Heckman's theorem \cite{DH}, $[\omega_t] = [\omega_0] - et$ where $e$ is the Euler class of $S^1$-fibration $\mu^{-1}(0) \rightarrow M_0$, where $M_0$ is the symplectic reduction at 0 with the induced symplectic form $\omega_0$. Hence we have

   \begin{equation}\label{Duist}
        \mathrm{DH}(t) = \int_{M_0} \frac{1}{(n-1)!}([\omega_0] - et)^{n-1}
   \end{equation}
   on $(c_1,c_2) \subset  \textrm{Im} \mu$.

   Note that a continuous function on an open interval $g : (a,b) \rightarrow \R$ is \textit{concave} if $g(tc + (1-t)d) \geq tg(c) + (1-t)g(d)$ for any $c,d \in (a,b)$ and for any $t \in (0,1)$. We remark the basic property of a concave function as follow.

   \begin{remark}\label{logconcave}
        Let $g$ be a continuous, piecewise smooth function on a connected interval $I \subset \R$. Then g is concave on $I$ if and only if
        the derivative of $g$ is decreasing, i.e. $g''(t) \leq 0$ for every smooth point $t \in I$ and $g'_+(c) - g'_-(c) < 0$ for every singular point $c \in I$, where $g'_+(c) = \lim_{t\rightarrow c, t > c} g'(t)$ and $g'_-(c) = \lim_{t\rightarrow c, t < c} g'(t)$.
   \end{remark}


   Note that Duistermaat and Heckman proved that $\mathrm{DH}$ is a polynomial on a connected regular open interval $U \subset \mu(M)$.
   The following formula due to Guillemin, Lerman and Sternberg describes the behavior of $\mathrm{DH}$ near the critical value of $\mu$.
   In particular, it implies that $\mathrm{DH}$ is $k$-times differentiable at a critical value $c \in \mu(M)$ if and only if
   $\mu^{-1}(c)$ does not contain a fixed component whose codimension is less than $4+2k$.

   \begin{theorem}\cite{GLS}\label{GLS}
       Assume that $c$ is a critical value which corresponds to the fixed components $C_{i}$'s. Then the jump of $DH(t)$ at $c$ is given by $$DH_+ - DH_- = \sum_{i} \frac{\textrm{vol}(C_i)}{(d-1)!\prod_ja_j} (t-c)^{d-1} + O((t-c)^d)$$
       where the sum is over the components $C_i$ of $M^{S^1} \bigcap \mu^{-1}(c)$, $d$ is half the real codimension of $C_i$ in $M$, and the $a_j$'s are the weights of the $S^1$-representation on the normal bundle of $C_i$.
   \end{theorem}
   If $c$ is a critical value which is not an extremum, then the codimension of the fixed point set in $\mu^{-1}(c)$ is at least 4. Therefore
   the theorem \ref{GLS} implies that $\mathrm{DH}(t)$ is continuous at non-extremal critical values and $\mathrm{DH}'(t)$ jumps at $c$ when $d$ equals 2. In the case
   when $d=2$, the two nonzero weights must have opposite signs, so the jump in the derivative is negative, i.e. $\mathrm{DH}'(t)$ decreases when it passes through the critical value with $d=2$. Since $\mathrm{DH}$ is continuous, the jump in $\frac{d}{dt}\ln{\mathrm{DH}(t)} = \frac{\mathrm{DH}'(t)}{\mathrm{DH}(t)}$ is negative at $c$. Combining with the lemma \ref{logconcave}, we have the following corollary.

   \begin{corollary}
       Let $(M,\omega)$ be a closed Hamiltonian $S^1$-manifold with the moment map $\mu : M \rightarrow \R$. Then the corresponding Duistermaat-Heckman function $\mathrm{DH}$ is log-concave on $\mu(M)$ if $(\log{\mathrm{DH}(t)})'' \leq 0$ for every regular value $t \in \mu(M)$.
   \end{corollary}

   Note that $(\log{\mathrm{DH}(t)})'' = \frac{\mathrm{DH}(t) \cdot \mathrm{DH}''(t) - \mathrm{DH}'(t)^2}{\mathrm{DH}(t)^2}$. Therefore
   $(\log{\mathrm{DH}(t)})'' \leq 0$ is equivalent to $\mathrm{DH}(t) \cdot \mathrm{DH}''(t) - \mathrm{DH}'(t)^2 \leq 0$.
   The equation (\ref{Duist}) implies that

   \begin{equation}\label{Duist2}
       \mathrm{DH}(t) \cdot \mathrm{DH}''(t) = (n-1)(n-2)\int_{M_0}e^2[\omega_t]^{n-3} \cdot \int_{M_0}[\omega_t]^{n-1}
   \end{equation}
   \begin{equation}\label{Duist3}
       \mathrm{DH}'(t)^2 = (n-1)^2\left(\int_{M_0}e[\omega_t]^{n-2}\right)^2
   \end{equation}

    To compute the integrals appeared in the equations (\ref{Duist2}) and (\ref{Duist3}), we need the following procedure :
    For an inclusion $\iota : \mu^{-1}(0) \hookrightarrow M$, we have a ring homomorphism $\kappa : H^*_{S^1}(M;\R) \rightarrow H^*_{S^1}(\mu^{-1}(0);\R) \cong H^*(M_0;\R)$ which is called the Kirwan map. Due to the Kirwan surjectivity \cite{Ki}, $\kappa$ is a ring surjection.

    Consider a 2-form $\widetilde{\omega} := \omega - d(\mu\theta)$ on $ES^1 \times M$ where $\theta$ is the connection form on $ES^1$. We denote by $x = \pi^*u \in H^2_{S^1}(M;\Z)$ where $\pi : M \times_{S^1} ES^1 \rightarrow BS^1$ and $u$ is a generator of $H^*(BS^1;\Z)$ such that
    the Euler class of the Hopf bundle $ES^1 \rightarrow BS^1$ is $-u$.
    Some part of the following two lemmas are given in \cite{Au}, but we give the complete proofs here.
    \begin{lemma}\label{equivariantform}
        $\widetilde{\omega}$ is $S^1$-invariant and closed, and $i_X \widetilde{\omega} = 0$ so that $\widetilde{\omega}$ represents a cohomology class in $H^*_{S^1}(M;\R)$. Moreover, for any fixed component $F \in M^{S^1}$, we have $\kappa([\widetilde{\omega}]) = [\omega_0]$ and $[\widetilde{\omega}]|_{F} = [\omega]|_F + \mu(F)u$. In particular, if $F$ is isolated, then $[\widetilde{\omega}]|_{F} = \mu(F)u$.
    \end{lemma}

    \begin{proof}
        For the first statement, it is enough to show that $i_X \widetilde{\omega}$ and $L_X \widetilde{\omega}$ vanish. Note that $i_X \widetilde{\omega} =
        i_X \omega - i_X d(\mu\theta) = -d\mu + di_X(\mu\theta) - L_X(\mu\theta)$ by Cartan's formula. Since $i_X(\mu\theta) = \mu$ and $\mu\theta$ is invariant under the circle action, we have $i_X \widetilde{\omega} = -d\mu + d\mu = 0$. Moreover, it is obvious that $\widetilde{\omega}$ is closed by definition. Hence $L_X \widetilde{\omega} = 0$ by Cartan's formula again.

        To prove the second statement, consider the following diagram.
        \begin{displaymath}
            \begin{array}{ccc}
                \mu^{-1}(0) \times ES^1 & \hookrightarrow & M \times ES^1\\
                \downarrow  &                 & \downarrow \\
                \mu^{-1}(0) \times_{S^1} ES^1 & \hookrightarrow & M \times_{S^1} ES^1 \\
                \downarrow_{\cong}  &                 & \\
                \mu^{-1}(0)/S^1 \cong M_{red}  &                 &  \\
            \end{array}
        \end{displaymath}
        Since $d\mu$ is zero on the tangent bundle $\mu^{-1}(0) \times ES^1$, the pull-back of $\widetilde{\omega} = \omega - d\mu \wedge \theta - \mu d\theta$ to $\mu^{-1}(0) \times ES^1$ is the restriction $\omega|_{\mu^{-1}(0) \times ES^1}$. And the lift-down of $\omega|_{\mu^{-1}(0) \times ES^1}$ to $\mu^{-1}(0)/S^1$ is just a reduced symplectic form at the level $0$. Hence $\kappa([\widetilde{\omega}]) = [\omega_0]$.

        To show the last statement, consider $[\widetilde{\omega}]|_{F} = [\omega - d(\mu\theta)]|_F = [\omega - d\mu\theta - \mu d\theta]|_F$. Since the restriction $d\mu|_{F \times ES^1}$ vanishes, we have $[\widetilde{\omega}]|_{F} = [\omega]|_{F} - \mu(F) \cdot [d\theta]|_{ES^1} = [\omega]|_{F} + \mu(F)u$. If $F$ is isolated, then we have $[\widetilde{\omega}]|_{F} = \mu(F)u$.

    \end{proof}

    \begin{lemma}\label{euler}
        Consider a 2-form $d\theta$ on $ES^1 \times M$ where $\theta$ is a connection 1-form on $ES^1$. Then we can lift $d\theta$ down to $ES^1 \times_{S^1} M$ so that $d\theta$ represents a cohomology class in $H^*_{S^1}(M;\R)$. Moreover, $[d\theta] = -x$ and $\kappa([d\theta]) = -\kappa(x) = e$ where $e$ is the Euler class of the $S^1$-fibration $\mu^{-1}(0) \rightarrow M_{red}$.
    \end{lemma}

    \begin{proof}
        Note that $i_X d\theta = L_X\theta - di_X\theta = 0$. Hence we can lift $d\theta$ down to $ES^1 \times_{S^1} M$. For any fixed point $x \in
        M^{S^1}$, the restriction $[d\theta]|_x$ is the Euler class of $ES^1 \times x \rightarrow BS^1$. Hence $[d\theta] = -u \cdot 1 = -x$. Here, the multiplication $\cdot $ comes from the $H^*(BS^1)$-module structure on $H^*_{S^1}(M)$. By the diagram in the proof of the Lemma \ref{equivariantform}, $\kappa([d\theta])$ is just the Euler class of the $S^1$-fibration $\mu^{-1}(0) \rightarrow \mu^{-1}(0)/S^1$. Therefore $\kappa([d\theta]) = -\kappa(x) = e$.
    \end{proof}

    Combining the equations (\ref{Duist2}), (\ref{Duist3}), Lemma \ref{equivariantform}, and Lemma \ref{euler}, we have the following corollary.
    \begin{corollary}\label{log2}
         $\mathrm{DH}(0)\cdot \mathrm{DH}''(0) - \mathrm{DH}'(0)^2 \leq 0$ if and only if $$(n-2)\int_{M_0} \kappa([d\theta]^2 [\widetilde{\omega}]^{n-3}) \cdot \int_{M_0}\kappa([\widetilde{\omega}]^{n-1}) - (n-1)(\int_{M_0}\kappa([d\theta] [\widetilde{\omega}]^{n-2}))^2 \leq 0.$$
    \end{corollary}

To compute the above integrals $\int_{M_0} \kappa([d\theta]^2 [\widetilde{\omega}]^{n-3}), \int_{M_0}\kappa([\widetilde{\omega}]^{n-1}),$ and $\int_{M_0}\kappa([d\theta] [\widetilde{\omega}]^{n-2})$, we need the residue formula due to Jeffrey and Kirwan. (See \cite{JK} and \cite{J}).

\begin{theorem}\cite{JK}\label{JK}
   Let $\nu \in H^*_{S^1}(M;\R)$. Then $$\int_{M_0} \kappa(\nu) = \sum_{F \in M^{S^1}, ~\mu(F)>0} Res (\frac{\nu|_F}{e_F}).$$
   Here $e_F$ is the equivariant Euler class of the normal bundle to F, and $Res(f)$ is a residue of $f$.
\end{theorem}

Now, let's compute $\int_{M_0} \kappa([d\theta]^2 [\widetilde{\omega}]^{n-3})$. By Theorem \ref{JK},
$$\int_{M_0} \kappa([d\theta]^2 [\widetilde{\omega}]^{n-3}) = \sum_{F \in M^{S^1}, ~\mu(F)>0} Res (\frac{[d\theta]^2 [\widetilde{\omega}]^{n-3}|_F}{e_F})$$

Since $[\widetilde{\omega}]|_z = \mu(z)u$ and $[d\theta]|_z = -u$ by lemma \ref{equivariantform} and \ref{euler}, we have
\begin{displaymath}
            \begin{array}{cl}
                \int_{M_0} \kappa([d\theta]^2 [\widetilde{\omega}]^{n-3}) & = \sum_{F \in M^{S^1}, ~\mu(F)>0} Res (\frac{\mu(F)^{n-3}u^{n-1}}{e_F})\\
                 & \\
                              &  =\sum_{F \in M^{S^1}, ~\mu(F)>0} Res(\frac{\mu(F)^{n-3}u^{n-1}}{m_Fu^n})\\
                 & \\
                              &  =\sum_{F \in M^{S^1}, ~\mu(F)>0} \frac{1}{m_F}\mu(F)^{n-3}\\
            \end{array}
\end{displaymath}
where $m_F$ is the product of all weights of tangential $S^1$-representation at $F$.
Similarly, if $\xi \in \R$ is a regular value of $\mu$, then we put $\widetilde{\mu} = \mu - \xi$ be the new moment map. By the same argument,
we have the following lemma

\begin{lemma}\label{eq}
    For a regular value $\xi$ of the moment map $\mu$,
    \begin{enumerate}
        \item
        \begin{displaymath}
        \int_{M_{\xi}} \kappa([d\theta]^2 [\widetilde{\omega}]^{n-3}) = \sum_{\substack{F \in M^{S^1},\\ \mu(F)>\xi}} \frac{1}{m_F}(\mu(F)-\xi)^{n-3}.
        \end{displaymath}
        \item
        \begin{displaymath}
        \int_{M_{\xi}} \kappa([d\theta] [\widetilde{\omega}]^{n-2}) = \sum_{\substack{F \in M^{S^1}, \\\mu(F)>\xi}} \frac{-1}{m_F}(\mu(F)-\xi)^{n-2}.
        \end{displaymath}
        \item
        \begin{displaymath}
        \int_{M_{\xi}} \kappa([\widetilde{\omega}]^{n-1}) = \sum_{\substack{F \in M^{S^1}, ~\\ \mu(F)>\xi}} \frac{1}{m_F}(\mu(F)-\xi)^{n-1}.
        \end{displaymath}
    \end{enumerate}
\end{lemma}

Combining the corollary \ref{log2} and the lemma \ref{eq}, we have the following proposition.

\begin{proposition}\label{det}
    Let $(M,\omega)$ be a closed Hamiltonian $S^1$-manifold with the moment map $\mu : M \rightarrow \R$. Assume that $M^{S^1}$ consists of isolated fixed points. Then a density function of the Duistermaat-Heckman measure with respect to $\mu$ is log-concave if and only if
    \begin{displaymath}
    \sum_{\substack{F \in M^{S^1}, \\ \mu(F)>\xi}} \frac{1}{m_F}(\mu(F)-\xi)^{n-3} \cdot \sum_{\substack{F \in M^{S^1}, \\ \mu(F)>\xi}} \frac{1}{m_F}(\mu(F)-\xi)^{n-1}
    \\ - \left(\sum_{\substack{F \in M^{S^1}, \\ \mu(F)>\xi}} \frac{1}{m_F}(\mu(F)-\xi)^{n-2}\right)^2 \leq 0
    \end{displaymath}
    for all regular value $\xi \in \mu(M)$, where $m_F$ is the product of all weights of the $S^1$-representation on $T_F M$. In particular, the log-concavity of the Duistermaat-Heckman measure is completely determined by
    the set $\{(\mu(F), m_F)_F | F \in M^{S^1} \}$.
\end{proposition}

\begin{corollary}\label{det2}
    Let $(M^{2n},\omega)$ and $(N^{2n}, \sigma)$ be two closed Hamiltonian $S^1$-manifold with the moment map $\mu_1$ and $\mu_2$ respectively.
    Assume there exist a bijection $\phi : M^{S^1} \rightarrow N^{S^1}$ which satisfies
    \begin{enumerate}
        \item for each $F \in M^{S^1}$, $m_F = m_{\phi(x)}$, and
        \item for each $F \in M^{S^1}$, $\mu_1(F) = \mu_2(\phi(F))$.
    \end{enumerate}
    where $m_F$ is the product of all weights of the tangential $S^1$-representation at $F$.
    If $N$ satisfies the log-concavity of the Duistermaat-Heckman measure with respect to $\mu_2$, then so does $M$ with respect to $\mu_1$.

\end{corollary}

\begin{remark}
    The integration formulae (1) and (3) in the lemma \ref{eq} are proved by Wu by using the stationary phase method. See Theorem 5.2 
    in \cite{Wu} for the detail. 
\end{remark}

\section{proof of the Theorem \ref{first}}

As noted in the introduction, if a Hamiltonian $S^1$-action on the K\"{a}hler manifold is holomorphic, then the corresponding Duistermaat-Heckman function is log-concave by \cite{Gr}. Let $(M^{2n},\omega)$ be a closed semifree Hamiltonian circle action with the moment map $\mu$. Assume that all fixed points are isolated. Let $\mathrm{DH}$ be the corresponding Duistermaat-Heckman function with respect to $\mu$. We will show that there is a K\"{a}hler form $\omega_1 \oplus \cdots \oplus \omega_n$ on $S^2 \times \cdots \times S^2$ with the standard diagonal holomorphic semifree circle action such that a bijection $\psi^{-1} \circ \phi : M^{S^1} \rightarrow (S^2 \times \cdots \times S^2)^{S^1} $ given in Section 2 satisfies the conditions in the corollary \ref{det2}, which implies the log-concavity of $\mathrm{DH}$. Now we start with the lemma below.

\begin{lemma}\label{determine}
   Let $(M^{2n},\omega)$ be a closed semifree Hamiltonian circle action with the moment map $\mu$. Assume that all fixed points are isolated. Then $\{(\mu(F), m_F)_F | F \in M^{S^1} \}$ is completely determined by $\mu(p_0^1), \mu(p_1^1), \cdots, \mu(p_1^n)$, where $p_k^j$'s are the fixed points of index $2k$ for $j=1,\cdots,{n \choose k}$.
\end{lemma}

\begin{proof}
    Consider an equivariant symplectic 2-form $\widetilde{\omega}$ on $ES^1 \times_{S^1} M$ which is given in Section 3. Because
    $x, a_1, \cdots, a_n$ form a basis of $H^2_{S^1}(M;\Z)$, we may let $[\widetilde{\omega}] = m_0x + m_1a_1 + \cdots + m_na_n$ for some real numbers $m_i$'s. (See Section 2 : the definition of $x, a_1, \cdots, a_n$.) By lemma \ref{equivariantform}, we have $[\widetilde{\omega}]|_{p_0^1} = \mu(p_0^1)u$. On the other hand, the right hand side is $(m_0x + m_1a_1 + \cdots + m_na_n)|_{p_0^1} = m_0u$, since every $a_i$ vanishes on $p_0^1$. Hence $m_0 = \mu(p_0^1)$. Similarly, $[\widetilde{\omega}]|_{p_1^i} = \mu(p_1^i)u$ and $(m_0x + m_1a_1 + \cdots + m_na_n)|_{p_1^i} = m_0u + m_iu$. Hence we have $m_i = \mu(p_1^i) - m_0 = \mu(p_1^i) - \mu(p_0^1)$ for each $i=1 , \cdots, n$. Therefore  $\mu(p_0^1), \mu(p_1^1), \cdots, \mu(p_1^n)$ determine the coefficients $m_i$'s of $[\widetilde{\omega}]$.
    For $p_k^j$, $[\widetilde{\omega}]|_{p_k^j} = \mu(p_k^j)u$ and $(m_0x + m_1a_1 + \cdots + m_na_n)|_{p_k^j} = m_0u + \sum_{i \in \phi(p_k^j)} m_iu$. Hence for fixed $k$, the set $\{(\mu(p_k^j),m_{p_k^j})_j | j=1,\cdots, {n \choose k} \}$ is just a $\{(m_0 + m_{i_1} + \cdots, + m_{i_k},(-1)^k)_{\{i_1, \cdots, i_k\}} | \{i_1, \cdots, i_k\} \subset \{1,2,\cdots,n\}\}$. Hence $\{(\mu(F),m_F)_F | F \in M^{S^1}\} = \cup_{k=0}^{k=n} \{(\mu(p_k^j),m_{p_k^j})_j | j=1,\cdots, {n \choose k} \}$ is completely determined by $\mu(p_0^1)$, $\mu(p_1^1)$, $\cdots, \mu(p_1^n)$.
\end{proof}

Now we are ready to prove the theorem \ref{first}.

\begin{proof}[\textit{Proof of Theorem \ref{first}}]
    For $\psi : (S^2 \times \cdots \times S^2)^{S^1} \rightarrow 2^{[n]}$ defined in Section 2, Let $\omega_i$ be the Fubini-Study form on $S^2$ such that the symplectic volume is $\mu(p_1^i) - \mu(p_0^1)$. Let $\mu' : S^2 \times \cdots \times S^2 \rightarrow \R$ be the moment map whose minimum is
    $\mu(p_0^1)$. Then we have $\{(\mu(F), m_F)_F | F \in M^{S^1} \} = \{(\mu(F), m_F)_F | F \in (S^2 \times \cdots \times S^2)^{S^1} \}$ and $\psi^{-1} \circ \phi : M^{S^1} \rightarrow (S^2 \times \cdots \times S^2)^{S^1}$ satisfies the condition in the corollary \ref{det2}.
    Therefore the Duistermaat-Heckman measure is log-concave on $\mu(M)$.

\end{proof}

\section{The hard Lefschetz property of the reduced symplectic forms}

For a K\"{a}hler manifold $(N,\sigma)$ with a holomorphic circle action preserving $\sigma$, its symplectic reduction is again K\"{a}hler, and hence
the reduced symplectic form $\sigma_t$ is K\"{a}hler. Therefore $\sigma_t$ satisfies the Hard Lefschetz property for every regular value $t$.
In this section, we show that the same thing happens when $(M,\omega)$ is a closed semifree Hamiltonian $S^1$-space whose fixed points are all isolated. The following theorem is due to Tolman and Weitsman.

\begin{theorem}\cite{TW2}\label{Kir}
    Let $S^1$ act on a compact symplectic manifold $(M,\omega)$ with moment map $\mu : M \rightarrow \R$. Assume that all fixed points are isolated and $0$ is a regular value. Let $M^{S^1}$ denote the set of fixed points.
    Define $K_+ := \{\alpha \in H^*_{S^1}(M;\Z) | \alpha_{F_+} = 0\}$ where $F_+ := M^{S^1} \bigcap \mu^{-1}(0,\infty)$ and $K_- := \{\alpha \in H^*_{S^1}(M;\Z) | \alpha_{F_-} = 0\}$ where $F_- := M^{S^1} \bigcap \mu^{-1}(-\infty,0)$.
    Then there is a short exact sequence:
    $$ 0 \longrightarrow K \longrightarrow H^*_{S^1}(M;\Z) \stackrel{\kappa}\longrightarrow (M_{red};\Z) \longrightarrow 0$$ where $\kappa : H^*_{S^1}(M;\Z) \rightarrow H^*(M_{red};\Z)$ is the Kirwan map.
\end{theorem}

\begin{proof}[\textit{Proof of Theorem \ref {second}}]
    Let $\kappa_M : H^*_{S^1}(M;\R) \rightarrow H^*(M_{red};\R)$ be the Kirwan map for $(M,\omega)$ and let $\kappa$ be the one for $(S^2 \times \cdots \times S^2, \sigma)$, where $\sigma := \omega_1 \oplus \cdots \oplus \omega_n$ is chosen in the proof of Theorem \ref{first} in Section 4.
    As in the proof of Theorem \ref{first}, we proved that there exists a semifree holomorphic Hamiltonian $S^1$-manifold $(S^2 \times \cdots \times S^2, \sigma)$ with the moment map $\mu'$ such that $\psi^{-1} \circ \phi : M^{S^1} \rightarrow (S^2 \times \cdots \times S^2)^{S^1}$ preserves their indices, weights, and the values of the moment map. Hence $\psi^{-1} \circ \phi$ identifies $K^M_+$ with $K^{S^2 \times \cdots \times S^2}_+$ and $K^M_-$ with $K^{S^2 \times \cdots \times S^2}_-$. The ring isomorphism $f : H^*_{S^1}(M;\Z) \rightarrow H^*_{S^1}(S^2 \times \cdots \times S^2;\Z)$ given in Corollary \ref{iso} satisfies  $\alpha|_x = f(\alpha)|_{\psi^{-1} \circ \phi(x)}$ for any $\alpha \in H^*_{S^1}(M;\Z)$ and any fixed point $x \in M^{S^1}$.
    Hence if $\alpha \in K^M_+$, then $f(\alpha) \in K^{S^2 \times \cdots \times S^2}_+$. Similarly for any $\alpha \in K^M_-$, we have $f(\alpha) \in K^{S^2 \times \cdots \times S^2}_-$. Therefore $f$ preserves the kernel of the Kirwan map $\kappa$ by Theorem \ref{Kir}.

    Note that $\kappa(f([\widetilde{\omega}]))$ is the reduced symplectic class of $S^2 \times \cdots \times S^2$ at 0. Since the K\"{a}hler quotient of the holomorphic action is again  K\"{a}hler, $\kappa(f([\widetilde{\omega}]))$ satisfies the hard Lefschetz property. Now, assume that $\omega_0$ does not satisfy the hard Lefschetz property. Then there exists a positive integer $k(<n)$ and some nonzero $\alpha \in H^k(M_{red};\R)$ such that $\alpha \cdot [\omega_0]^{n-k} = 0$ in $H^{2n-k}(M_{red})$. By the Kirwan surjectivity, we can find $\widetilde{\alpha} \in H^*_{S^1}(M;\R)$
    with $\kappa(\widetilde{\alpha}) = \alpha$. Then $\widetilde{\alpha} \cdot [\widetilde{\omega}]^{n-k}$ is in $\ker{\kappa}$ and hence the image
    $f(\widetilde{\alpha} \cdot [\widetilde{\omega}]^{n-k})$ is in $\ker{\kappa}$, since $f$ maps $\ker{\kappa}$ of $M$ to $\ker{\kappa}$ of $S^2 \times \cdots \times S^2$. It implies that $f(\widetilde{\alpha}) = 0$ by the hard Lefschetz condition for $f([\widetilde{\omega}])$. It contradicts that $f$ is an isomorphism.

    It remains to show that $(M,\omega)$ satisfies the Hard Lefschetz property. Remind that $\psi^{-1} \circ \phi : M^{S^1} \rightarrow (S^2 \times \cdots \times S^2)^{S^1}$ induces an isomorphism $$f : H^*_{S^1}(M;\R) \rightarrow H^*_{S^1}(S^2 \times \cdots \times S^2;\R),$$ which sends the equivariant symplectic class $[\widetilde{\omega}]$ to $[\widetilde{\sigma}]$ as we have seen in Section 4. Here, $\widetilde{\sigma}$ is an equivariant symplectic form induced by $\sigma - d(\mu'\theta)$ in $H^*_{S^1}(S^2 \times \cdots \times S^2;\R)$.
    Since $f$ is a $H^*(BS^1;\R)$-isomorphism, it induces a ring isomorphism
    \begin{displaymath}
        f_u : \frac{H^*_{S^1}(M;\R)}{u \cdot H^*_{S^1}(M;\R)} \rightarrow \frac{H^*_{S^1}(S^2 \times \cdots \times S^2;\R)}{u \cdot H^*_{S^1}(S^2 \times \cdots \times S^2;\R)}
    \end{displaymath}
    Moreover, the quotient map $\pi_M : H^*_{S^1}(M;\R) \rightarrow \frac{H^*_{S^1}(M;\R)}{u \cdot H^*_{S^1}(M;\R)} \cong H^*(M;\R)$ ($\pi_{S^2 \times \cdots \times S^2}$ respectively) is just a ring homomorphism which comes from an inclusion $M \hookrightarrow M \times_{S^1} ES^1$ as a fiber. Therefore $\pi_M([\widetilde{\omega}]) = [\omega]$ and $\pi_{S^2 \times \cdots \times S^2}([\widetilde{\sigma}]) = [\sigma]$. It means that
    $f_u : H^*(M;\R) \rightarrow H^*(S^2 \times \cdots \times S^2 ;\R)$ is a ring isomorphism which sends $[\omega]$ to $[\sigma]$. Since $\sigma$ is a K\"{a}hler form, it satisfies the Hard Lefschetz property. Hence so $\omega$ does.

\end{proof}

\bigskip
\bibliographystyle{amsalpha}

\end{document}